\documentclass{conm-p-l}
\usepackage{amsmath,amssymb}
\usepackage{graphicx}
\usepackage{verbatim}
\usepackage[square]{natbib}
\usepackage{enumerate}

\numberwithin{equation}{section}

\newtheorem{thm}{Theorem}[section]
\newtheorem{lem}[thm]{Lemma}
\newtheorem{pr}[thm]{Proposition}
\newtheorem{cor}[thm]{Corollary}
\newtheorem{defn}[thm]{Definition}

\newtheorem*{unremark}{Remark}

\newcommand{\Cox}{\hfill \Box}

\newcommand{\Em}[1]{\textbf{#1}}

\newcommand{\ee}{\epsilon}

\newcommand{\vv}{\mathbf{v}}
\newcommand{\xx}{\mathbf{x}}
\newcommand{\yy}{\mathbf{y}}
\newcommand{\zz}{\mathbf{z}}
\newcommand{\pp}{\mathbf{p}}

\newcommand{\ww}{\mathbf{w}}
\newcommand{\rr}{\mathbf{r}}

\newcommand{\I}{\mathcal{I}}
\newcommand{\CC}{\mathcal{C}}

\newcommand{\DD}{\mathcal{D}}
\newcommand{\M}{\mathcal{M}}
\newcommand{\J}{\mathbf{J}}
\newcommand{\one}{\mathbf{1}}
\newcommand{\zero}{\mathbf{0}}

\newcommand{\C}{\mathbb C}

\newcommand{\R}{\mathbb R}

\newcommand{\denom}{H}

\newcommand{\nbd}{\mathcal{N}}

\DeclareMathOperator{\grad}{\nabla}
\DeclareMathOperator{\Real}{Re}
\DeclareMathOperator{\Imag}{Im}

\DeclareMathOperator{\hess}{\mathcal{H}}

\begin{document}

\title{Asymptotic expansions of oscillatory integrals with 
complex phase}
\author{Robin Pemantle}
\address{University of Pennsylvania, Department of Mathematics, 
209 S. 33rd Street, Philadelphia, PA 19104 USA}
\email{pemantle@math.upenn.edu}
\thanks{Research supported in part by National Science Foundation grant 
\# DMS 0603821}

\author{Mark C. Wilson}
\address{Department of Computer Science, University of Auckland, 
Private Bag 92019, Auckland, NEW ZEALAND}
\email{mcw@cs.auckland.ac.nz} 



\begin{abstract}
We consider saddle point integrals in $d$ variables whose
phase function is neither real nor purely imaginary.  
Results analogous to those for Laplace (real phase) and 
Fourier (imaginary phase) integrals hold whenever the
phase function is analytic and nondegenerate.  These results
generalize what is well known for integrals of Laplace and Fourier
type.  The method is via contour shifting in complex $d$-space.
This work is motivated by applications to asymptotic enumeration.
\end{abstract}

\maketitle

\section{Introduction} \label{sec:intro}

Integrals of the form
\begin{equation} \label{eq:I}
\I (\lambda) := \I (\lambda ; \phi , A) := 
   \int e^{-\lambda \phi (x)} A (x) \, dx
\end{equation}
arise in many areas of mathematics.  There are many variations.
This integral can arise in one or more variables; the variables
may be real or complex; the integral may be global or taken in
over a small neighborhood or oddly shaped set; varying degrees
of smoothness may be assumed; varying degrees of degeneracy may
be allowed near the critical points of the phase function, $\phi$.
Often what is sought is a leading order estimate of $\I (\lambda)$
as the positive real parameter $\lambda$ tends to $\infty$, or
an asymptotic series $\I (\lambda) \sim \sum_n c_n g_n (\lambda)$
where $\{ g_n \}$ is a given sequence of elementary functions 
and the expansion is possibly nowhere convergent but satisfies
\begin{equation} \label{eq:expansion}
\I (\lambda) - \sum_{n=0}^{N-1} c_n g_n (\lambda) = O(g_N (\lambda))
\end{equation}
for any $N \geq 1$.  

In recent work on the asymptotics of multivariate generating
functions~\cite{PW1,PW2,PW9,BBBP,BP-cones,BP-residues},
we have required results of this type in the case where the
phase function $\phi$ and the amplitude function $A$ are 
analytic functions of several variables.  The phase function
is typically neither real nor purely imaginary, but satisfies
$\Real \{ \phi \} \geq 0$.  Although these results are known,
at least at a folklore level, we could not find this case
analyzed in the literature in the generality we required.
We therefore wrote this note to patch the gap and to remedy
what would have been ghost citations in~\cite[Lemmas~4.7~and~4.8]{PW2}.

To better explain the relation between the present results and
the literature, we discuss two special cases in which the results
are well known, but by substantially different methods.

\subsection*{Integrals of Fourier type}

Let $f, A : \R^d \to \R$ be smooth (that is, $C^\infty$) functions 
of $d$ real variables, with $A$ having compact support.  
Taking $\phi = -i f$ gives the Fourier-type integral
$$\I (\lambda) = \int e^{i \lambda f(x)} A(x) \, dx \, .$$
The standard method of studying this integral is as follows.
If $f$ has no critical points in the support of $A$, then 
integration by parts shows $\I (\lambda)$ to be rapidly decreasing:
$\I (\lambda) = O(\lambda^{-N})$ for all positive integers, $N$.
Using a partition of unity, the integral may therefore be localized
to neighborhoods of the critical points of $f$.  At an isolated
critical point $\grad f$ vanishes; if the Hessian matrix is non-degenerate
then the Morse lemma produces a smooth change of variables
under which $f(x) = S_\pm (x) := \sum_{j=1}^d \pm x_j^2$.  This reduces 
the general problem to the case $\phi = i S$.  To solve this, 
expand $A$ in a series.  Each term may be explicitly integrated,
resulting in an expansion in decreasing powers of $\lambda$:
\begin{equation} \label{eq:fourier asym}
\I (\lambda) = \sum_{n \geq 0} c_n \, \lambda^{-d/2 - n} \, .
\end{equation}
To see that the resulting series for $\I (\lambda)$ 
satisfies~\eqref{eq:expansion}, one uses integration by parts
again to bound the remainder term.  The coefficients $\{ c_n \}$ 
are determined by the derivatives of $A$ and $\phi$.  A particularly 
lucid treatment of this may be found in the first two sections 
of~\cite[Chapter~VIII]{stein}.  

\subsection*{Integrals of Laplace type}

The other well studied case is the Laplace-type integral, where
$\phi$ is real.  Localization of the integral to the minima
of $\phi$ is immediate because the integrand is exponentially
small elsewhere.  Integrating over balls whose radius has order 
$\lambda^{-1/2}$ shows that near a quadratically nondegenerate 
minimum, $x_0$, the bound $A(x) = O(|x - x_0|^N)$ translates into
the bound 
\begin{equation} \label{eq:O}
\I (\lambda) = O(\lambda^{-(d+N)/2}) \, .
\end{equation}
Again, one changes variables, expands $A$ into a power series,
and integrates term by term, to obtain the series
\begin{equation} \label{eq:laplace asym}
\I (\lambda) = \sum_{n \geq 0} c_n \, \lambda^{-(d+n)/2} \, .
\end{equation}
Applying~\eqref{eq:O} to the integral of the remainder term 
shows that~\eqref{eq:laplace asym} is an asymptotic expansion 
for $\I (\lambda)$.  Classical treatments of integrals of Laplace 
type may be found in many 
places~\cite{bleistein-handelsman,wong-asymptotic-integrals},
often accompanied by separate treatments of the Fourier case.

\subsection*{Complex methods}

The series~\eqref{eq:fourier asym} and~\eqref{eq:laplace asym}
are formally identical.  Further inspection of formulae for $c_n$ 
in the literature shows these to be nearly the same as well, 
differing only by constant factors of unit modulus.  This points 
to the possibility of unifying the results and generalizing to 
arbitrary complex functions.  An assumption of analyticity
is required; technically, this is stronger than smoothness but
in practice one is never satisfied without the other.  Assuming 
analyticity, derivations in the Fourier and Laplace cases may 
indeed by unified via a hybrid approach.  In one variable, this 
is carried out all the time under various names such as 
``steepest descent'', ``stationary phase'' or ``saddle point''.  
Writing $\I (\lambda)$ as a complex contour integral, 
the critical point will be a saddle point in $\C^1$ 
for the real part of the phase function; the contour may then 
be re-oriented to pass through the saddle in the direction of 
steepest descent of $-\phi$, converting the integral into 
one of Fourier type and explaining why the series are nearly
identical.  This is carried out, for example, in~\cite{debruijn} 
or~\cite[Chapter~7]{bleistein-handelsman}.  

In more than one variable, complex phase results 
generalizing~\eqref{eq:fourier asym} and~\eqref{eq:laplace asym}
are strangely absent from the literature.  They are not
treated in~\cite{stein,henrici-2,bleistein-handelsman,
wong-asymptotic-integrals,breitung,debruijn}.  We have found 
only two significant treatments of complex phases.  One is 
Section~7.7 of~\cite{hormander-partial-differential1}.  Here, 
notions from several complex variable theory are used to validate
integration by parts, leading to asymptotic expansions in the case 
when $A$ is smooth and compactly supported.  This work has been
used by~\cite{raichev-wilson-higher-order} to produce 
explicit higher order terms in various expansions.  As stated,
H\"ormander's method does not allow integration over manifolds
with boundary.  The second treatment is~\cite{fedoryuk}, which 
uses similar methodology to ours (complex contour deformation).
The cited work omits the details, referring to the original 
article~\cite{fedoryuk-book} in Russian.  Regions of integration 
with boundaries are allowed here, but only if $\Real \{ \phi \}$ 
is strictly maximized on the interior (our preliminary case, 
discussed in Section~\ref{sec:strict}).  In applications, work 
requiring high-dimensional saddle integration with complex phases 
is usually done by hand; see, e.g.,~\cite{wormald-tournaments} 
or~\cite{BFSS-short}.

There are several possible reasons for this gap.  One is that 
fewer people are aware of the corresponding techniques in 
several complex variables.  Results such as Cauchy's integral
theorem in several variables and Stokes' theorem for 
holomorphic $d$-forms in $\C^d$, while known to first-year
graduate students in several complex variables, are not always
known to those applying these techniques.  The few mathematicians
who do use these techniques to evaluate integrals~\cite{howls,lichtin} 
are not unusually interested in anything as pedestrian as the 
generalization of standard saddle point integral evaluation.
Secondly, even if one knows the techniques, moving contours 
in higher dimensions is a tricky business.  The geometry can
be much more complicated.  When the real part of the phase
vanishes at points that are not critical, further care must be 
taken in how the contour is moved.
Finally, the need for such a generalization, illuminating 
though it may be, may have been too infrequent.

The organization of rest of the paper is as follows.  In the next 
section we give some definitions having to do with stratified 
spaces.  We then state our main result, Theorem~\ref{th:main}. 
Section~\ref{sec:preliminaries} records some easy computations
in the case $\phi$ is the standard phase function $S(\xx) := 
\sum_{j=1}^d x_j^2$.  Section~\ref{sec:strict} handles
the general case under the assumption that the real part
of $\phi$ has a strict minimum at the critical point.  
Theorem~\ref{th:main} is proved in Section~\ref{sec:proofs}.  
Section~\ref{sec:examples} gives an application from~\cite{PW2}
to the estimation of coefficients of a bivariate generating function.
Finally, Section~\ref{sec:comments} discusses further research
directions.

\setcounter{equation}{0}
\section{Notation and statement of results} \label{sec:notation}

\subsection*{Stratified spaces}

Because of the useful properties ensuing from the definition, we will 
use Whitney stratified spaces as our chains of integration.  Aside
from these useful properties, the details of the definition need 
not concern us, though for completeness we give a precise definition.
Let $I$ be a finite partially ordered set and define an $I$-decomposition 
of a topological space $Z$ to be a partition of $Z$ into a disjoint 
union of sets $\{ S_\alpha : \alpha \in I \}$ such that
$$S_\alpha \cap \overline{S_\beta} \neq \emptyset
          \iff S_\alpha \subseteq \overline{S_\beta}
          \iff \alpha \leq \beta \, .$$
\begin{defn}[Whitney stratification] \label{def:whitney 2}
Let $Z$ be a closed subset of a smooth manifold $\M$.
A \Em{Whitney stratification} of $Z$ is an $I$-decomposition such that
\begin{enumerate} [(i)]
\item Each $S_\alpha$ is a manifold in $\R^n$.
\item If $\alpha < \beta$, if the sequences $\{ x_i  \in S_\beta \}$
and $\{ y_i \in S_\alpha \}$ both converge to $y \in S_\alpha$, if the
lines $l_i = \overline{x_i \, y_i}$ converge to a line $l$ and the
tangent planes $T_{x_i} (S_\beta)$ converge to a plane $T$ of some
dimension, then both $l$ and $T_y (S_\alpha)$ are contained in $T$.
\end{enumerate}
\end{defn}
For example, any manifold is a Whitney stratified space with one
stratum; any manifold with boundary is a Whitney stratified space
with two strata, one being the interior and one the boundary; a
$k$-simplex is a Whitney stratified space whose strata are all
its faces.  Whitney stratified spaces are closed under 
products and the set of products of strata will stratify the
product.  Every algebraic variety admits a Whitney stratification,
although the singular locus filtration may be too coarse to be
a Whitney stratification.

\subsection*{Critical points}

Associated with the definition of a stratification is the
stratified notion of a critical point.  Observe that under 
this definition, any zero-dimensional stratum of $\M$ is 
a critical point of $\M$.  
\begin{defn}[smooth functions and their critical points]
\label{def:crit}
Say that a function $\phi : \M \to \C$ on a stratified space $\M$
is smooth if it is smooth when restricted to each stratum.  A point
$p \in \M$ is said to be critical for the smooth function $\phi$
if and only if the restriction $d\phi_{|S}$ vanishes, where $S$ is
the stratum containing $p$.  
\end{defn}
Let $\M \subseteq \C^d$ be a real analytic, $d$-dimensional
stratified space.
This means that each stratum $S$ is a subset of $\C^d$ and 
each of the chart maps $\psi$ from a neighborhood of the 
origin in $\R^k$ to some $k$-dimensional stratum 
$S \subseteq \C^d$ is analytic (the coordinate functions 
are convergent power series) with a nonsingular differential.  
It follows that $\psi$ may be extended to a holomorphic map 
on a neighborhood of the origin in $\C^k$, whose range we denote 
by $S \otimes \C$.  Choosing a small enough neighborhood, we may 
arrange for $S \otimes \C$ to be a complex $k$-manifold embedded 
in $\C^d$.  We say a function $f : \M \to \C$ is analytic if
it has a convergent power series expansion in a neighborhood
of every point; an analytic function on $\M$ may be extended
to a complex analytic function on a neighborhood of $\M$ in
$\M \otimes \C := \bigcup_{\alpha \in I} S_\alpha \otimes \C$.
Because we are interested in the integrals of $d$-forms over $\M$, 
there is no loss of generality in assuming that $\M$ is contained
in the closure of its $d$-dimensional strata, whence $\M \otimes \C$
is a neighborhood of $\M$ in $\C^d$.  Real $d$-manifolds in $\C^d$ are
not naturally oriented, so we must assume that an orientation is
given for each $d$-stratum of $\M$, meaning that the chart maps
from $\R^d$ to $\M$ must preserve orientations.
The critical point $p \in \M$ is said to be \Em{quadratically nondegenerate} 
if $p$ is in a $d$-dimensional stratum and 
$\hess (p)$ is nonsingular, where 
$$\hess (p) := \left ( \frac{\partial^2 \phi}{\partial x_i \, \partial x_j}
   \right )_{1 \leq i,j \leq d}$$
is the Hessian matrix for the function $\phi$ of Definition~\ref{def:crit}
on a neighborhood of $p$ in $\C^d$.

Finally, we define a point of $\M$ to be \Em{stationary} if it is 
critical and $\Real \{ \phi \}$ vanishes there.

\subsection*{Results}

The main result of the paper is an asymptotic expansion for $\I (\lambda)$.
\begin{thm} \label{th:main}
Let $\M$ be a compact, real analytic stratified space of dimension $d$ 
embedded in $\C^d$.  Let $A$ and $\phi$ be analytic functions on 
a neighborhood of $\M$ and suppose $\Real \{ \phi \} \geq 0$ on $\M$.  
Let $G$ be the subset of stationary points of $\phi$ on $M$ and assume that 
$G$ is finite and each stationary point is quadratically nondegenerate.  Then 
the integral $\I (\lambda) := \int_\M e^{-\lambda \phi (x)} A(x) \, dx$ has
an asymptotic expansion
$$\I (\lambda) \sim \sum_{\ell = 0}^\infty c_\ell \lambda^{-(d+\ell)/2} 
   \, .$$
If $A$ is nonzero at some point of $G$ then the leading term is
given by
\begin{equation} \label{eq:leading}
c_0 = (2 \pi)^{-d/2} \sum_{x \in G} A(x) \, \omega (x , \lambda)
   \left ( \det \hess (x) \right )^{-1/2} \, .
\end{equation}
where $\omega (x , \lambda)$ is the unit complex number 
$e^{- \lambda \phi (x)}$.  The choice of sign of the square root
on the right-hand side of~\eqref{eq:leading} is determined by 
choosing any analytic chart map $\psi$ for a neighborhood of $x$
and defining
$$(\det \hess (x))^{-1/2} := \left ( \det \hess (\phi \circ \psi) 
   \right )^{-1/2} \, \det J_\psi$$
where $J_\psi$ is the Jacobian matrix for $\psi$ at $x$ and 
the $-1/2$ power of the determinant on the right is the product 
of the inverses of principal square roots of the eigenvalues.
\end{thm}

\begin{cor} \label{cor:half}
Assume the hypotheses of Theorem~\ref{th:main} except that some 
of the points of $G$ may be in strata of dimension $d-1$ with
a neighborhood in $\M$ diffeomorphic to a halfspace in $\R^d$.
Then the same conclusion holds except that the summand 
in~\eqref{eq:leading} corresponding to such a point $x$ must be 
multiplied by $1/2$.
\end{cor}



\setcounter{equation}{0}
\section{Preliminary results for the standard phase function} 
\label{sec:preliminaries}

For $\xx \in \R^d$, let $S(\xx) := \sum_{j=1}^d x_j^2$ 
denote the standard quadratic form.   We begin with a couple 
of results on integrals of Laplace type where the phase 
function is the standard quadratic and $A$ is a monomial.
In one dimension, 
\begin{pr} \label{pr:lambda}
$$\int_{-\infty}^\infty x^{n} \, e^{-\lambda x^2} \, dx =
   \frac{n!}{(n/2)! \, 2^n \, \sqrt{2 \pi}} \, \lambda^{-1/2 - n/2}$$
if $n$ is even, while the integral is zero if $n$ is odd.
\end{pr}

\begin{proof} For odd $n$ the result is obvious from the fact
that the integrand is an odd function.  For $n=0$ the result
is just the standard Gaussian integral.  By induction, assume 
now the result for $n-2$.  Integrating by parts to get the 
second line, we have 
\begin{eqnarray*}
\int x^{n} e^{- \lambda x^2} \, dx & = &
   \int \frac{-x^{n-1}}{2 \lambda} 
   \left ( -2 \lambda x \, e^{-x^2} \, dx \right )   \\
& = & \frac{n-1}{2 \lambda} \int x^{n-2} \, e^{-x^2} \, dx \\
& = & \frac{n-1}{2 \lambda} \frac{1}{\sqrt{2 \pi}}
   \frac{(n-2)!}{(n/2-1)! \, 2^{n-1}} \lambda^{1/2 - n/2)}
\end{eqnarray*}
by the induction hypothesis.  This is equal to
$\lambda^{-1/2 - n/2} (2 \pi)^{-1/2} n! \, / \, ((n/2)! \, 2^{n})$, 
completing the induction.
\end{proof}

\begin{cor}[monomial integral] \label{cor:mv monomial}
Let $\rr$ be any $d$-vector of nonnegative integers and let
$\xx^\rr$ denote the monomial $x_1^{r_1} \cdots x_d^{r_d}$.  Then
$$\int_{\R^d} \xx^{\rr} e^{- \lambda S(\xx)} \, d\xx =
  \beta_\rr \lambda^{-(d + |\rr|) / 2}$$
where 
\begin{equation} \label{eq:beta}
\beta_\rr := (2 \pi)^{-d/2} 
   \prod_{j=1}^d \frac{r_j !}{(r_j /2)! \, 2^{r_j}}
\end{equation}
if all the components $r_j$ are even, and zero otherwise.
\end{cor}

\begin{proof} The integral factors into
$$\prod_{j=1}^d \left [  \int_{-\infty}^\infty x_j^{r_j} \,
   e^{- \lambda r_j^2} \, dx_j \right ] \, .$$
\end{proof}

To integrate term by term over a Taylor series for $A$, we
need the following estimate.

\begin{lem}[big-O estimate] \label{lem:big-O}
Let $A$ be any smooth function satisfying $A(\xx) = O(|\xx|^r)$
at the origin.  Then the integral of $A(\xx) e^{-\lambda S(\xx)}$
over any compact set $K$ may be bounded from above as
$$\int_K A(\xx) e^{- \lambda S(\xx)} \, d\xx = O(\lambda^{-(d+r)/2})$$
\end{lem}

\begin{proof} Because $K$ is compact and $A(\xx) = O(|\xx|^r)$
at the origin, it follows that there is some constant $c$ for which
$|A(\xx)| \leq c |\xx|^r$ on all of $K$.
Let $K_0$ denote the intersection of $K$ with the ball
$|\xx| \leq \lambda^{-1/2}$ and for $n \geq 1$ let $K_n$ denote
the intersection of $K$ with the shell $2^{n-1} \lambda^{-1/2}
\leq |\xx| \leq 2^n \lambda^{-1/2}$.  On $K_0$ we have
$$|A(\xx)| \leq c \lambda^{-r/2}$$
while trivially
$$\int_{K_0} e^{-\lambda S(\xx)} \, d\xx \leq \int_{K_0} d\xx
   \leq c_d \lambda^{-d/2} \, .$$
Thus
$$\left | \int_{K_0} A(\xx) e^{- \lambda S(\xx)} \, d\xx \right |
   \leq c' \lambda^{-(r+d)/2} \, .$$
For $n \geq 1$, on $K \cap K_n$, we have the upper bounds
\begin{eqnarray*}
|A(\xx)| & \leq 2^{rn} & c \lambda^{-r/2} \\[1ex]
e^{- \lambda S(\xx)} & \leq & e^{-2^{n-1}} \\[1ex]
\int_{K_n} \, d\xx & \leq & 2^{rn} \, c_d \, \lambda^{-d/2} \, .
\end{eqnarray*}
Letting $c'' := c \cdot c_d \cdot \sum_{n=1}^\infty 2^{2rn} \,
e^{-2^{n-1}} < \infty$, we may sum to find that
$$\sum_{n=0}^\infty \left | \int_{K_n} A(\xx)
   e^{-\lambda S(\xx)} \, d\xx \right | \leq (c' + c'')
   \lambda^{-(r+d)/2} \, ,$$
proving the lemma.
\end{proof}

It is now easy to compute a series for $\I (\lambda)$ in
the case where $\phi$ is the standard quadratic and the
integral is over a neighborhood of the origin in $\R^d$.
\begin{thm}[standard phase] \label{th:standard}
Let $A (\xx) = \sum_\rr a_\rr \xx^\rr$ be any analytic function
defined on a neighborhood $\nbd$ of the origin in $\R^d$.  Let
\begin{equation} \label{eq:I lam}
\I (\lambda) := \int_\nbd A(\xx) e^{-\lambda S(\xx)} \, d\xx \, .
\end{equation}
Then
$$\I (\lambda) \sim \sum_n \left[\sum_{|\rr| = n} a_\rr \beta_\rr\right]
   \lambda^{-(n + d)/2}$$
as an asymptotic series expansion in decreasing powers of $\lambda$,
with $\beta_\rr$ as in~\eqref{eq:beta}.
\end{thm}

\begin{proof} Write $A(\xx)$ as a power series up to 
degree $N$ plus a remainder term:
$$A(\xx) = \left ( \sum_{n=0}^N \sum_{|\rr| = n} a_\rr \xx^\rr
   \right ) + R(\xx)$$
where $R(\xx) = O(|\xx|^{N+1})$.  Using Corollary~\ref{cor:mv monomial}
to integrate all the monomial terms and Lemma~\ref{lem:big-O} to
bound the integral of $R(\xx) e^{-\lambda S(\xx)}$ shows that
$$\I (\lambda) = \sum_{n=0}^N \sum_{|\rr| = n} a_\rr \beta_\rr
   \lambda^{-(d+n)/2} + O(\lambda^{-(d+n)/2 - 1})$$
which proves the asymptotic expansion.
\end{proof}

\setcounter{equation}{0}
\section{The case of a strict minimum} \label{sec:strict}

In this section, we continue to integrate over a neighborhood
of the origin in $\R^d$, but we generalize to any analytic
phase function $\phi$ with the restriction that the real part
of $\phi$ have a strict minimum at the origin.  The assumption
of a strict minimum localizes the integral to the origin, so
the only tricky aspects are keeping track of the sign
(Lemma~\ref{lem:sign}) and being rigorous about moving
the contour.
\begin{thm} \label{th:strict}
Let $A$ and $\phi$ be complex-valued analytic functions on a compact 
neighborhood $\nbd$ of the origin in $\R^d$ and suppose that the 
real part of $\phi$ is nonnegative, vanishing only at the origin.
Suppose that the Hessian matrix $\hess$ of $\phi$ at the origin
is nonsingular.  Denoting $\I (\lambda) := \int_\nbd A(\xx) 
e^{-\lambda \phi (\xx)} \, d\xx$, there is an asymptotic expansion
$$\I (\lambda) \sim \sum_{\ell \geq 0} c_\ell \lambda^{d/2 - \ell}$$
where
$$c_0 = A(\zero) \frac{(2 \pi)^{-d/2}}{\sqrt{\det \hess}}$$
and the choice of sign is defined by taking the product of
the principal square roots of the eigenvalues of $\hess$.
\end{thm}

The proof is essentially a reduction to the case of standard 
phase.  The key is the well known Morse Lemma.  The proof
given in~\cite[VIII:2.3.2]{stein} is for the smooth category
and purely real or imaginary phase but extends without
significant change to complex values and the analytic
category.  For completeness, we include the adapted proof.

\begin{lem}[Morse lemma] \label{lem:morse}
There is a bi-holomorphic change of variables $\xx = \psi (\yy)$
such that $\phi (\psi (\yy)) = S(\yy) := \sum_{j=1}^d y_j^2$.
The differential $J_\psi = d\psi (\zero)$ will satisfy $(\det J_\psi)^2 =
(\det \frac{1}{2} \hess )^{-1}$.
\end{lem} 

\begin{proof} Addressing the second conclusion first, we recall 
how the Hessian matrix behaves under a change of variables.  
If $\psi : \C^d \to \C^d$ is bi-holomorphic on a neighborhood 
of $\xx$ and if $\phi$ has vanishing gradient at $\psi (\xx)$ 
and Hessian matrix $\hess$ there, then the Hessian matrix 
$\tilde{\hess}$ of $\phi \circ \psi$ at $\xx$ is given by
$$\tilde{\hess} = J_\psi^T \hess J_\psi$$
where $J_\psi$ is the Jacobian matrix for $\psi$ at $\xx$.  The
standard form $S$ has Hessian matrix equal to twice the identity,
hence any function $\psi$ satisfying $\phi \circ \psi = S$ must satisfy 
$$2 \, Id = J_\psi^T \hess J_\psi \, .$$
Dividing by two and taking determinants yields 
$|\J_\psi|^2 \det ( \frac{1}{2} \hess ) = 1$, proving the second
conclusion.

To prove the change of variables, the first step is to write
$$\phi (\xx) = \sum_{j,k = 1}^d x_j x_k \phi_{j,k}$$
where the functions $\phi_{j,k} = \phi_{k,j}$ are analytic and satisfy
$\phi_{j,k} (\zero) = (1/2) \hess_{j,k}$.  It is obvious from a
formal power series viewpoint that this can be done because
the summand $x_j x_k \phi_{j,k}$ can be any power series with
coefficients indexed by the orthant $\{ \rr : \rr \geq \delta_j +
\delta_k \}$; these orthants cover $\{ \rr : |\rr| \geq 2 \}$,
so we may obtain any function $\phi$ vanishing to order two;
matching coefficients on the terms of order precisely two shows
that $\phi_{j,k} (\zero) = (1/2) \hess_{j,k}$.

More constructively, we may give a formula for $\phi_{j,k}$.
There is plenty of freedom, but a convenient choice is to let
$a_\rr$ denote the coefficient of $\xx^\rr$ in $\phi (\xx)$ and 
to take
$$x_k x_k \phi_{j,k} (\xx) := \sum_{|\rr| \geq 2}
   \frac{r_j (r_k - \delta_{j,k})}{|\rr| (|\rr| - 1)} \; a_\rr \xx^\rr \, .$$
For fixed $\rr$, it is easy to check that
$$\sum_{1 \leq j , k \leq d} \frac{r_j (r_k - \delta_{j,k})}
   {|\rr| (|\rr| - 1)} = 1$$
whence $\phi = \sum x_j x_k \phi_{j,k}$.  Alternatively, the
following analytic computation from~\cite{stein} verifies that
$\phi = \sum_{j,k} x_j x_k \phi_{j,k}$.  Any function $f$
vanishing at zero satisfies $f(t) = \int_0^1 (1-s) f' (s) \, ds$,
as may be seen by integrating by parts (take $g(s) = - (1-s)$).
Fix $\xx$ and apply this with $f(t) = (d/dt) \phi (t \xx)$ to obtain
$$\phi (\xx) = \int_0^1 \frac{d}{dt} \phi (t \xx) \, dt =
   \int_0^1 (1-t) \frac{d^2}{dt^2} \phi (t \xx) \, dt \, .$$
The multivariate chain rule gives
$$\frac{d^2}{dt^2} \phi (t \xx) = \sum_{j,k} x_j x_k
   \frac{\partial^2 \phi}{\partial x_j \partial x_k} (t \xx) \, ;$$
plug in $\phi = \sum_\rr a_\rr \xx^\rr$ and integrate term by term
using $\int_0^1 (1-t) t^{n-2} \, dt = \frac{1}{n(n-1)}$ to see
that $\phi = \sum_{j,k} x_j x_k \phi_{j,k}$.

The second step is an induction.  Suppose first that $\phi_{j,j} (\zero)
\neq 0$ for all $j$.  The function $\phi_{1,1}^{-1}$ and a branch of
the function $\phi_{1,1}^{1/2}$ are analytic in a neighborhood of
the origin.  Set
$$y_1 := \phi_{1,1}^{1/2} \left [ x_1 + \sum_{k > 1}
   \frac{y_k \phi_{1,k}}{\phi_{1,1}} \right ] \; .$$
Expanding, we find that the terms of $y_1^2$ of total degree
at most one in the terms $x_2 , \ldots , x_d$ match those of
$\phi$ and therefore,
\begin{equation} \label{eq:y_1}
\phi (\xx) = y_1^2 + \sum_{j,k \geq 2} x_j x_k h_{j,k}
\end{equation}
for some analytic functions $h_{j,k}$ satisfying $h_{j,k} (\zero)
= (1/2) \hess_{j,k}$.  Similarly, if
$$\phi (\xx) = \sum_{j=1}^{r-1} y_j^2 + \sum_{j,k \geq r} x_j x_k h_{j,k}$$
then setting
$$y_r := \phi_{r,r}^{1/2} \left [ x_r + \sum_{k > r}
   \frac{y_k \phi_{r,k}}{\phi_{r,r}} \right ]$$
gives
$$\phi (\xx) = \sum_{j=1}^{r} y_j^2
   + \sum_{j,k \geq r+1} x_j x_k \tilde{h}_{j,k}$$
for some analytic functions $\tilde{h}_{j,k}$ still satisfying
$h_{j,k} (\zero) = (1/2) \hess_{j,k}$.  By induction, we arrive
at $\phi (\xx) = \sum_{j=1}^d y_j^2$, finishing the proof of
the Morse Lemma in the case where each $\hess_{j,j}$ is nonzero.

Finally, if some $\hess_{j,k} = 0$, because $\hess$ is nonsingular
we may always find some unitary map $U$ such that the Hessian
$U^T \hess U$ of $\phi \circ U$ has no vanishing diagonal entries.
We know there is a $\psi_0$ such that $(\phi \circ U) \circ \psi_0 = S$,
and taking $\psi = U \circ \psi_0$ finishes the proof in this case.
\end{proof}

\noindent{\sc Proof of Theorem}~\ref{th:strict}: 
The power series allows us to extend $\phi$
to a neighborhood of the origin in $\C^d$.  Under the change
of variables $\psi$ from the previous lemma, we see that
\begin{eqnarray*}
\I (\lambda) & = & \int_{\psi^{-1} \CC} A \circ \psi (\yy) e^{-\lambda S(\yy)}
   (\det d\psi (\yy)) \, d\yy \\
& := & \int_{\psi^{-1} \CC} \tilde{A} (\yy) e^{-\lambda S(\yy)} \, d\yy
\end{eqnarray*}
for some analytic function $\tilde{A}$, 
where $\CC$ is a neighborhood of the origin in $\R^n$.  We need to 
check that we can move the chain $\psi^{-1} \CC$ of integration back 
to the real plane.

Let $h(\zz) := \Real \{ S(\zz) \}$.  The chain $\CC' := \psi^{-1} (\CC)$
lies in the region $\{ \zz \in \C^d : h(\zz) > 0 \}$
except when $\zz = 0$, and in particular, $h \geq \ee > 0$
on $\partial \CC'$.  Let
$$H(\zz , t) := \Real \{ \zz \} + (1-t) \, i \, \Imag \{ \zz \} \, .$$
In other words, $H$ is a homotopy from the identity map to the
map $\pi$ projecting out the imaginary part of the vector $\zz$.
For any chain $\sigma$, the homotopy $H$ induces a chain homotopy,
$H(\sigma)$ supported on the image of the support of $\sigma$
under the homotopy $H$ and satisfying
$$\partial H(\sigma) = \sigma - \pi \sigma + H(\partial \sigma) \, .$$
With $\sigma = \CC'$, observing that $S(H(\zz , t)) \geq S(\zz)$,
we see there is a $(d+1)$-chain $\DD$ with
$$\partial \DD = \CC' - \pi \CC' + \CC''$$
and $\CC''$ supported on $\{ h > \ee \}$.  Stokes Theorem tells us
that for any holomorphic $d$-form $\omega$,
$$\int_{\partial \DD} \omega = \int_\DD d\omega = 0$$
and consequently, that
$$\int_{\CC'} \omega = \int_{\pi \CC'} \omega + \int_{\CC''} \omega \, .$$
When $\omega = \tilde{A} e^{-\lambda S} \, d\yy$, the integral over
$\CC''$ is $O(e^{-\lambda \ee})$, giving
$$\I (\lambda) = \int_{\pi \CC'} \tilde{A} (\yy) e^{-\lambda S(\yy)} \, d\yy
   + O(e^{-\ee \lambda}) \, .$$

Up to sign, the chain $\pi \CC''$ is a disk in $\R^d$ with the
standard orientation plus something supported in $\{ h > \ee \}$.
To see this, note that $\pi$ maps any real $d$-manifold in $\C^d$
diffeomorphically to $\R^d$ wherever the tangent space is transverse
to the imaginary subspace.  The tangent space to the support of $\CC'$
at the origin is transverse to the imaginary subspace because $S \geq 0$
on $\CC'$, whereas the imaginary subspace is precisely the negative
$d$-space of the index-$d$ form $S$.  The tangent space varies
continuously, so in a neighborhood of the origin, $\pi$ is a
diffeomorphism.  Observing that $\tilde{A} (\zero) = A(\zero)
\det (d\psi(\zero)) = A(\zero) (\det \hess)^{-1/2}$ and using
Theorem~\ref{th:standard} finishes the proof up to the choice of 
sign of the square root.

The map $d\pi \circ d\psi^{-1} (\zero)$ maps the standard
basis of $\R^d$ to another basis for $\R^d$.  Verifying the
sign choice is equivalent to showing that this second basis
is positively oriented if and only if $\det (d\psi (\zero))$
is the product of the principal square roots of the eigenvalues
of $\hess$ (it must be either this or its negative).  Thus we will
be finished by applying the following lemma (with $\alpha = \psi^{-1}$).
\begin{lem} \label{lem:sign}
Let $W \subseteq \C^d$ be the set $\{ \zz : \Real \{ S(\zz) \} > 0 \}$.
Pick any $\alpha \in GL_d (\C)$ mapping $\R^d$ into $\overline{W}$
and let $M := \alpha^T \alpha$ be the matrix representing
$S \circ \alpha$.  Let $\pi : \C^d \to \R^d$ be projection
onto the real part.  Then $\pi \circ \alpha$ is orientation
preserving on $\R^d$ if and only if $\det \alpha$ is the
product of the principal square roots of the eigenvalues
of $M$ (rather than the negative of this).
\end{lem}

\begin{proof} First suppose $\alpha \in GL_d (\R)$.
Then $M$ has positive eigenvalues, so the product of their
principal square roots is positive.  The map $\pi$ is the
identity on $\R^d$ so the statement boils down to saying
that $\alpha$ preserves orientation if and only if it
has positive determinant, which is true by definition.
In the general case, let $\alpha_t := \pi_t \circ \alpha$,
where $\pi_t (\zz) = \Real \{ \zz \} + (1-t) \Imag \{ \zz \}$.
As we saw in the previous proof, $\pi_t (\R^d) \subseteq \overline{W}$
for all $0 \leq t \leq 1$, whence $M_t := \alpha_t^T \alpha_t$
has eigenvalues with nonnegative real parts.  The product of the
principal square roots of the eigenvalues is a continuous
function on the set of nonsingular matrices with no negative
real eigenvalues.  The determinant of $\alpha_t$ is a continuous
function of $t$, and we have seen it agrees with the product
of principal square roots of eigenvalues of $M_t$ when $t=1$
(the real case), so by continuity, this is the correct sign
choice for all $0 \leq t \leq 1$; taking $t=0$ proves the lemma.
\end{proof}

For later use, we record one easy corollary of Theorem~\ref{th:strict}.
\begin{cor} \label{cor:strict half}
Assume the hypotheses of Theorem~\ref{th:strict} and let
$\nbd'$ be the intersection of $\nbd$ with a region diffeomorphic
to a halfspace through the origin.  If $A(\zero) \neq 0$ then
$$\I' (\lambda) := \int_{\nbd'} A(\xx) e^{-\lambda \phi (\xx)} 
   \, dx \sim \frac{c_0}{2} \, \lambda^{-d/2}$$
where $c_0$ is the same as in the conclusion of Theorem~\ref{th:strict}.
\end{cor}

\begin{proof} Under the change of variables $\psi$
and the projection $\pi$, this region maps to a region $\nbd''$
diffeomorphic to a halfspace with the origin on the boundary.
Changing variables by $\yy = \lambda^{-1/2} \xx$ and
writing $\nbd_\lambda$ for $\lambda^{1/2} \nbd''$, we have
$$\I' (\lambda) = \lambda^{-d/2} \int_{\nbd_\lambda} A_\lambda (\yy)
   e^{-S(\yy)} \, d\yy$$
where $A_\lambda (\yy) = (A \circ \psi) (\lambda^{-1/2} \yy)$.
The function $A_\lambda$ converges to $A(\zero)$ pointwise
but also in $L^2 (\mu)$ where $\mu$ is the Gaussian measure
$e^{-S(\xx)} \, d\xx$.  Also, the regions $\nbd_\lambda$ converge 
to a halfspace $H$ in the sense that their indicators 
$\one_{\nbd_\lambda}$ converge to $\one_H$ in $L^2 (\mu)$.
Thus $A_\lambda \one_{\nbd_\lambda}$ converges to 
$A(\zero) \one_H$ in $L^1 (\mu)$, and unravelling this statement
we see that
$$\int_{\nbd_\lambda} A_\lambda (\yy) e^{-S(\yy)} \, d\yy
   \to \int_H A(\zero) e^{-S(\yy)}\, d\yy \, .$$
The last quantity is equal to $c_0 / 2$, showing that
$\lambda^{d/2} \I' (\lambda) \to c_0 / 2$ and finishing the proof.
\end{proof}

\setcounter{equation}{0}
\section{Proofs of main results} \label{sec:proofs}

Theorem~\ref{th:main} differs from Theorem~\ref{th:strict}
in several ways.  The most important is that the set where
$\Re \phi$ vanishes may extend to the boundary
of the region of integration.  This precludes the use of
the easy deformation $\pi$ because $\CC''$ is no longer
supported on $\{ h > \ee \}$.  Consequently, some work
is required to construct a suitable deformation.  We do so
via notions from stratified Morse theory~\cite{GM}.

\subsection*{Tangent vector fields}

If $\xx$ is a point of the stratum $S$ of the stratified space
$\M$, let $T_\xx (\M)$ denote the tangent space to $S$ at $\xx$.
Because $\M$ is embedded in $\C^d$, the tangent spaces may
all be identified as subspaces of $\C^d$.  Thus we have a notion
of the tangent bundle $T \M$, a section of which is simply
a vector field $f$ on $\M \subseteq \C^d$ such that 
$f (\xx) \in T_\xx (\M)$ for all $\xx$.  A consequence of
the two Whitney conditions is the local product structure 
of a stratified space: a point $\pp$ in a $k$-dimensional
stratum $S$ of a stratified space $\M$ has a neighborhood in which
$\M$ is homeomorphic to some product $S \times X$.  
According to~\cite{GM}, a proof may be found in mimoegraphed 
notes of Mather from 1970; it is based on Thom's Isotopy Lemma 
which takes up fifty pages of the same mimeographed notes.
The next lemma is the only place where we use this (or any)
consequence of Whitney stratification. 
\begin{lem} \label{lem:section}
Let $f$ be a smooth section of the tangent bundle to $S$, that is
$f(s) \in T_s (S)$ for $s \in S$.  Then each $s \in S$ has
a neighborhood in $\M$ on which $f$ may be extended to a smooth
section of the tangent bundle. 
\end{lem}

\begin{proof}
Parametrize $\M$ locally by $S \times X$ and 
extend $f$ by $f(s,x) := f(s)$.
\end{proof}

\begin{lem}[vector field near a non-critical point]
\label{lem:local V}
Let $\xx$ be a point of the stratum $S$ of the stratified space
$\M$ and suppose $\xx$ is not critical for the function $\phi$.
Then there is a vector $\vv \in T_\xx (S \otimes \M)$ such that
$\Real \{ d\phi (\vv) \} > 0$ at $\xx$.  Furthermore, there is a
continuous section $f$ of the tangent bundle in a neighborhood
$\nbd$ of $\xx$ such that $\Real \{ d\phi (f(\yy)) \} > 0$
at every $\yy \in \nbd$.
\end{lem}

\begin{proof} By non-criticality of $\xx$, there is a
$\ww \in T_\xx (S)$ with $d\phi (\ww) = u \neq 0$ at $\xx$.
Multiply $\ww$ componentwise by $\overline{u}$ to obtain $\vv$
with $\Real \{ d\phi (\vv) \} > 0$ at $\xx$.  Use any chart map
for $S \otimes \C$ near $\xx$ to give a locally trivial coordinatization
for the tangent bundle and define a section $f$ to be the constant
vector $\vv$; then $\Real \{ d\phi (f (\yy)) \} > 0$ on some
sufficiently small neighborhood of $\xx$ in $S$.  Finally, extend
to a neighborhood of $\xx$ in $\M$ by Lemma~\ref{lem:section}.
\end{proof}

Although we are working in the analytic category, the chains
of integration are topological objects, for which we may use
$C^\infty$ methods (in what follows, even $C^1$ methods will do).
In particular, a partition of unity argument enhances the local
result above to a global result.

\begin{lem}[global vector field, in the absence of critical points]
\label{lem:V_0}
Let $\M$ be a compact stratified space and $\phi$ a smooth
function on $\M$ with no critical points.  Then there is a global
section $f$ of the tangent bundle of $\M$ such that the real part
of $d\phi (f)$ is everywhere positive.
\end{lem}

\begin{proof} For each point $\xx \in \M$, let $f_\xx$
be a section as in the conclusion of Lemma~\ref{lem:local V},
on a neighborhood $U_\xx$.  Cover the compact space $\M$ by finitely
many sets $\{ U_\xx : \xx \in F \}$ and let $\{ \psi_\xx : \xx \in F \}$
be a smooth partition of unity subordinate to this finite cover.  Define
$$f(\yy) = \sum_{\xx \in F} \psi_\xx (\yy) f_\xx (\yy) \, .$$
Then $f$ is smooth; it is a section of the tangent bundle because
each tangent space is linearly closed; the real part of
$d\phi (f(\yy))$ is positive because we took a convex combination
in which each contribution was nonnegative and at least one was
positive.
\end{proof}

Another partition argument gives our final version of this result.
\begin{lem}[global vector field, vanishing only at critical points]
\label{lem:V}
Let $\M$ be a compact stratified space and $\phi$ a smooth function
on $\M$ with finitely many critical points.  Then there is a global
section $f$ of the tangent bundle of $\M$ such that the real part
of $d\phi (f)$ is nonnegative and vanishes only when $\yy$ is
a critical point.
\end{lem}

\begin{proof} Let $\M_\ee$ be the compact stratified space
resulting in the removal of an $\ee$-ball around each critical
point of $\phi$.  Let $f_\ee$ be a vector field as in the conclusion
of Lemma~\ref{lem:V} with $\M$ replaced by $\M_\ee$.  Let $c_n$
be a positive real number, small enough so that the magnitudes
of all partial derivatives of $c_n f_{1/n}$ of order up to $n$
are at most $2^{-n}$.  In the topology of uniform convergence of
derivatives of bounded order, the series $\sum_n c_n f_n$ converges
to a vector field $f$ with the required properties.
\end{proof}

\subsection*{Proof of Theorem~\protect{\ref{th:main}}}

Let $f$ be a tangent vector field along which $\phi$ increases
away from critical points, as given by Lemma~\ref{lem:V}.  
Such a field gives rise to a differential flow, which, informally, 
is the solution to $d \pp / dt = f(\pp)$.  To be more formal, 
let $\xx$ be a point in a stratum $S$ of $\M$.
Via a chart map in a neighborhood of $\xx$, we solve the ODE
$d \Phi(t) / dt = f(\Phi (t))$ with initial condition
$\Phi (0) = \xx$, obtaining a trajectory $\Phi$ on some interval
$[0,\ee_\xx]$ that is supported on $S$.  Doing this simultaneously
for all $\xx \in \M$ results in a map
$$\Phi : \M \times [0,\ee] \to \C^d$$
with $\Phi (\xx , t)$ remaining in $S \otimes \C$ when $\xx$ is in
the stratum $S$.  The map $\Phi$ satisfies $\Phi (\xx , 0) = \xx$
and $(d/dt) \Phi (\xx , t) = f (\Phi (\xx , t))$.  The fact that
this may be defined up to time $\ee$ for some $\ee > 0$ is a
consequence of the fact that the vector field $f$ is bounded
and that a small neighborhood of $\M$ in $\M \otimes \C$ is embedded
in $\C^d$.  Because $f$ is smooth and bounded, for sufficiently
small $\ee$ the map $\xx \mapsto \Phi (\xx , \ee)$ is a 
diffeomorphism.   

The flow reduces the real part of $\phi$ everywhere except the
critical points which are rest points.  Consequently, it defines
a homotopy $H(\xx , t) := \Phi (\xx , \ee t)$ between $\CC$ and a
chain $\CC'$ on which the minima of the real part of $\phi$ occur
precisely on the set $G$.  Recall that $H$ induces a chain homotopy
$\CC_H$ with $\partial \CC_H = \CC' - \CC + \partial \CC \times \sigma$,
where $\sigma$ is a standard 1-simplex.  Let $\omega$ denote the
holomorphic $d$-form $A(\zz) \exp (- \lambda \phi (\zz)) \, d\zz$.
Because $\omega$ is a holomorphic $d$-form in $\C^d$, we have
$d\omega = 0$.  Now, by Stokes' Theorem,
\begin{eqnarray*}
0 & = & \int_{\CC_H} d\omega \\
& = & \int_{\partial \CC_H} \omega \\
& = & \int_{\CC'} \omega - \int_{\CC} \omega -
   \int_{\partial \CC \times \sigma}  \omega .
\end{eqnarray*}
The chain $\partial \CC \times \sigma$ is supported on a finite union
of spaces $S \otimes \CC$ where $S$ is a stratum of dimension
at most $d-1$.  The integral of a holomorphic $d$-form
vanishes over such a chain.  Therefore, the last term
on the right drops out and we have
$$\int_{\CC} \omega = \int_{\CC'} \omega \, .$$

Outside of a neighborhood of $G$ the magnitude of the integrand
is exponentially small, so we have shown that there are $d$-chains
$\CC_\xx$ supported on arbitrarily small neighborhoods $\nbd (\xx)$
of each $\xx \in G$ such that
\begin{equation} \label{eq:localized}
\I (\lambda) - \sum_{\xx \in G} \int_{\CC_\xx} \omega
\end{equation}
is exponentially small.  To finish that proof, we need only
show that each $\int_{\CC_\xx} \omega$ has an asymptotic
series in decreasing powers of $\lambda$
whose leading term, when $A(\xx) \neq 0$, is given by
\begin{equation} \label{eq:leading 2}
c_0 (\xx) = (2 \pi)^{-d/2} A(\xx) e^{\lambda \phi (\xx)}
   (\det \hess (\xx))^{-1/2} \, .
\end{equation}
The $d$-chain $\CC_\xx$ may by parametrized by a map 
$\psi_\xx : B \to \nbd (\xx)$, mapping the origin to $\xx$,
where $B$ is the open unit ball in $\R^d$.  By the chain rule,
$$\int_{\CC_\xx} \omega = \int_B [A \circ \psi] (\xx)  
   \exp (- \lambda [\phi \circ \psi (\xx)]) \det d\psi (\xx) \, d\xx .$$
The real part of the analytic phase function $\phi \circ \psi$ has 
a strict minimum at the origin, so we may apply Theorem~\ref{th:strict}.
We obtain an asymptotic expansion whose first term is
\begin{equation} \label{eq:hess 2}
(2 \pi \lambda)^{-d/2} [A \circ \psi] (\zero) \det d\psi (\zero) (\det M_\xx)^{-1/2}
\end{equation}
where $M_\xx$ is the Hessian matrix of the function $\phi \circ \psi$.
The term $[A \circ \psi] (\zero)$ is equal to $A(\xx)$.
The Hessian matrix of $\phi \circ \psi$ at the origin is given by
$M_\xx = d\psi (\zero) \, \hess (\xx) \, d\psi (\zero)$.  Thus
$$\det M_\xx = (\det d\psi (\zero))^2 \det \hess (\xx)$$
and plugging into~\eqref{eq:hess 2} yields~\eqref{eq:leading 2}.
$\Cox$ \\[1ex]

\noindent{Proof of Corollary~\protect{\ref{cor:half}}}:
Lemma~\ref{lem:V} does not require the critical points
to be in the interior, so the argument leading up 
to~\eqref{eq:localized} is still valid.  For those 
points $\xx$ in a $(d-1)$-dimensional stratum, use 
Corollary~\ref{cor:half} in place of Theorem~\ref{th:strict}
to obtain~\eqref{eq:leading 2} with an extra factor of $1/2$.
$\Cox$

\begin{unremark}
The reason we do not continue with a litany of special
geometries (quarter-spaces, octants, and so forth) is that
the case of a halfspace is somewhat special.  The differential
of the change of variables at the origin is a nonsingular map,
which must send half-spaces to half-spaces, though it will
in general alter angles of any smaller cone.
\end{unremark}

\setcounter{equation}{0}
\section{Examples} \label{sec:examples}

The simplest multidimensional application of our results is 
a computation from~\cite{PW2}.  The purpose is to estimate
coefficients of a class of bivariate generating functions
whose denominator is the product of two smooth divisors.
We give only a brief summary of how one arrives at~\eqref{eq:PW2}
from a problem involving generating functions; a complete
explanation of this can be found in~\cite[Section~4]{PW2}.  
Note, however, that the mathematics of the integral is not contained 
in that paper, which instead refers to an earlier draft of this one!  

Let $v_1, v_2$ be distinct analytic functions of $z$ with $v_1(1) 
= v_2(1)=1, 0\neq v_1'(1)\neq v_2'(1)\neq 0$ and such that each $|v_i(z)|$ 
attains its maximum on $|z| = 1$ only at $z = 1$. 
For example, the last condition is satisfied by any pair of aperiodic power series with 
nonnegative coefficients and radius of convergence greater than $1$.

Consider the generating function $F(z,w) = 1/\denom(z,w)$, where 
$\denom(z,w) = (1-wv_1(z))(1-wv_2(z))$. The two branches of the 
curve $\denom = 0$ intersect only at $(1,1)$ and this intersection 
is transverse.  The Maclaurin coefficients
of $F(z,w) = \sum_{r,s} a_{rs} z^r w^s$ are given by the 
Cauchy integral formula
$$ a_{rs} = \frac{1}{(2\pi i)^2} \frac{dwdz}
   {z^{r+1}w^{s+1}(1-wv_1(z))(1-wv_2(z))} $$ 
where the integral is taken over a product of circles centred at $(0,0)$ 
and of sufficiently small radii.

Pushing the contour out to $|z| = 1, |w| = 1 - \varepsilon$ we 
obtain the same formula, since $F$ is still analytic inside the 
product of disks bounded by these latter circles.  Pushing the 
$w$-contour out to $|w| = 1 + \varepsilon$, using the residue formula 
on the inner integral and observing that the integral over 
$|w| = 1+ \varepsilon$ is exponentially decaying as $s\to \infty$, 
we see that
$$a_{rs} \approx \frac{1}{2\pi}\int_{|z| = 1} \frac{-R_s(z)}{z^{r+1}} 
   \, dz $$
where $\approx$ means that the difference is exponentially decaying 
as $s \to \infty$ and $R_s(z)$ denotes the sum of residues of 
$w \mapsto w^{-(s+1)} F(z,w)$ at the roots $w= 1/v_i(z)$, $i \in \{1, 2\}$.

The residue sum $R_s(z)$ can be rewritten in terms of an integral via 
$$-R_s(z) = (s+1)\int_0^1 \left[(1-p)v_1 + pv_2\right]^s \, dp$$
and so we have 
$$
a_{rs} \approx \frac{s+1}{2\pi} 
\int_{|z| = 1} z^{-(r+1)} \int_0^1 \left[(1-p)v_1(z) + pv_2(z)\right]^s \, dp dz.
$$
In order to cast this into our standard framework, we need to be able to define 
a branch of the logarithm of $(1-p)v_1(z) + pv_2(z)$. We do this by 
localizing on the circle $|z|=1$ to a sufficiently small 
neighbourhood of the point $z=1$. This is possible since the integrand decays 
exponentially away from $z=1$, by hypotheses on the $v_i$, and we shall show 
that the integral near $z=1$ decays only polynomially.

The substitution $z = e^{i\theta}$ converts this to an integral
\begin{equation} \label{eq:PW2}
a_{rs} \approx \frac{s+1}{2\pi} \int_{\nbd} \int_0^1 
   e^{-s \phi(p,t)} A(p,z) \, dp dt
\end{equation}
where $\phi(p,t) = ir\theta/s + \log \left[(1-p)v_1(e^{it}) 
+ pv_2(e^{it})\right]$, $A(p,t) = 1$, and $\nbd$ is a closed 
interval centred at $0$.  To compute asymptotics in the direction 
$r/s = \kappa$, for fixed $\kappa > 0$, we can consider $\phi$ to 
be independent of $r$ and $s$.  

We now asymptotically evaluate~\eqref{eq:PW2} using Theorem~\ref{th:main}.
We can rewrite the iterated integral as a single integral over the 
stratified space $\M = \nbd \times [0,1]$. The phase $\phi$ has 
nonnegative real part and this fits into our framework.  There 
is a single stationary point, at $(p,z) = (1/2, 0)$ 
(note that $\Real \{ \phi \}$ is zero for all $(p, 0)$, so 
Theorem~\ref{th:strict} does not suffice).  This critical point 
is quadratically nondegenerate and direct computation using 
Theorem~\ref{th:main} yields 
\begin{equation} \label{eq:c}
a_{rs} = \frac{1}{|v_1'(1) - v_2'(1)|} + O(s^{-1})
\end{equation}
as $s\to \infty$ with $\kappa$ fixed. By keeping track of error terms more 
explicitly, it is easily shown that this approximation is 
uniform in $\kappa$ provided $\kappa$ stays in a compact subset 
of the open interval formed by $v_1'(1), v_2'(1)$ (it follows from 
our assumptions that these numbers are positive real --- see \cite{PW2} 
for more details).  This  means that $a_{rs}$ is asymptotically constant 
in any compact subcone of directions away from the boundary formed by 
the lines $\kappa_i = v_i'(1)$.

This example, and in fact a number of cases in~\cite{PW2}, can also
be solved using iterated residues.  This is carried out
in~\cite{BP-residues}.  Iterated residues have the advantage
of showing that the $O(s^{-1})$ term decays exponentially,
but the disadvantage that they do not give any results when
$\kappa$ approaches the boundary.  The present methods do
give boundary results.  Corollary~\ref{cor:half} shows that
$a_{rs}$ converges to one half the right-hand side of~\eqref{eq:c}
when $(r,s) \to \infty$ with $r/s = \kappa_1 + O(1)$, and a small
extension obtains a Gaussian limit: letting $\Phi$ denote the
standard normal cumulative distribution function, we have
$$a_{rs} = \frac{\Phi (u)}{|v_1'(1) - v_2'(1)|} + O(s^{-1})$$
when $r,s \to \infty$ with $(r/s-\kappa_1) / s^{1/2} \to u$.

\section{Further comments}\label{sec:comments}

We have not emphasized explicit formulae for the higher order terms,
giving an equation such as~\eqref{eq:leading} only for the leading
term in the case where $A(\zero) \neq 0$.  However, our results 
establish the validity of existing computations of higher order terms 
under our more general hypotheses.  

To elaborate, we prove Theorem~\ref{th:main} by first constructing 
a change of variables $\xx \mapsto \Phi (x,\ee)$ homotopic to the 
identity under which the
minimum of $\Real \{ \phi \}$ at $\zero$ is strict, and then changing
variables, again homotopically to the identity, to the standard form. 
The composition $\psi$ of these two maps is homotopic to the identity
but is far from explicitly given: while the second map is constructed
by an explicitly defined Morse function, the first deformation is the
solution to a differential equation and is not particularly explicit.

In \cite{hormander}, H\"{o}rmander derives such an explicit formula 
(assuming smoothness) for integrals of our type where $\M = \R^d$ 
and $A$ has compact support.  The formula is indeed rewritten and 
used in \cite{raichev-wilson-higher-order} to compute higher order 
terms for generating function applications, in which more restrictive
hypotheses preclude the vanishing of $\Real \{ \phi \}$ on a curve
reaching the boundary of the chain of integration.  Their methods,
while not covering the cases of interest here, do have the virtue 
of dealing with the change of variables $\psi$ only through the
equation $S = \phi \circ \psi$.  In particular, the derivatives
of $\psi$ arising in the computation of the new amplitude function
$(A \circ \psi) \det d\psi$ can be computed by implicitly differentiating
the equation $S = \phi \circ \psi$.  Having found at least one such 
$\psi$ homotopic to the identity, we are now free to replicate 
the computations of~\cite{raichev-wilson-higher-order} under our 
more general hypotheses, as follows.

In the case of standard phase, the coefficient of $\lambda^{-(n+d)/2}$ 
is given by (provided all $r_i$ are even) 
$$
\sum_{|\rr| = n} a_\rr \beta_\rr
$$
where $a_\rr$ is the Maclaurin coefficient of $A$ corresponding 
to the monomial $\rr$ and $\beta_\rr$ is the constant defined in 
Corollary~\ref{cor:mv monomial}.  Note that $n$ must be even for 
this coefficient to be nonzero, so we write 
$n=2k$. The differential operator $\partial_1^{r_1} \dots \partial_d^{r_d}$ 
when applied to $A$ and evaluated at $\zero$ yields precisely 
$\prod_i r_i! a_\rr$. Thus the operator 
$$ \sum_{|\rr| = k} \frac{\partial_1^{2r_1} \dots \partial_d^{2r_d}}
   {4^k r_1! \dots r_d!} $$
applied to $A$ and evaluated at zero yields the coefficient we seek.

After the Morse lemma is applied using the change of variables 
$S=\phi\circ \psi$, we need to apply the displayed operator to 
the new amplitude $(A\circ \psi) \det d\psi$.  The resulting 
expression evaluated at $\xx$ can be computed directly via the rules of
Leibniz and Faa di Bruno.  Evaluating at $\xx$ simplifies some terms,
and, as mentioned above, derivatives of $(A \circ \psi) \det d\psi$ 
may be computed without explicitly specifying $\psi$.

As a relatively simple example, consider the case $k=1$ and $d=1$.
The differential operator reduces to $\frac{1}{4} \partial^2$ 
where $\partial$ denotes differentiation with respect to the variable $x$. 
Applying this to $(A\circ \psi) \det d\psi$ yields (with superscripts 
denoting the order of derivatives and arguments suppressed)
$$ \frac{1}{4} \left(A^{(2)} \left(\psi^{(1)}\right)^3 
   + 3A^{(1)}\psi^{(1)}\psi^{(2)} + A^{(0)} \psi^{(3)}\right).  $$
The defining equation $S = \phi \circ \psi$ can be differentiated 
to yield the system
\begin{align*}
2x & = \phi^{(1)} \psi^{(1)} \\
2 & = \phi^{(2)} \left[\psi^{(1)}\right]^2 + \phi^{(1)}\psi^{(2)}\\
0 & = \phi^{(3)} \left[\psi^{(1)}\right]^3 + 3\phi^{(2)}\psi^{(1)}\psi^{(2)} + \phi^{(1)} \psi^{(3)}\\
0 & = \phi^{(4)} \left[\psi^{(1)}\right]^4 + 6\phi^{(3)}\left[\psi^{(1)}\right]^2\psi^{(2)} 
+ 4 \phi^{(2)}\psi^{(1)}\psi^{(3)} + 3 \phi^{(2)}\left[\psi^{(2)}\right]^2 + \phi^{(1)} \psi^{(4)}.
\end{align*}

Evaluating these at the point in question, we see that the terms 
with highest derivatives of $\psi$ vanish in each equation. 
The system is triangular and can be solved explicitly to obtain
\begin{align*}
\psi^{(1)} & = \sqrt{\frac{2}{\phi^{(2)}}} \\
\psi^{(2)} & = \frac{- 2 \phi^{(3)}}{3 \left[\phi^{(2)}\right]^2}\\
\psi^{(3)} & = \frac{\left[5 \phi^{(3)}\right]^2 - 3\phi^{(2)}\phi^{(4)}}{3\sqrt{2}\left[\phi^{(2)}\right]^{7/2}}. 
\end{align*}
Putting these together with the expression for the derivative of 
$(A\circ \psi) \det d\psi$ above yields an expression for the 
$\lambda^{-3/2}$ term in the integral that is a rational function 
with denominator $\left[\phi^{(2)}\right]^{7/2}$ and numerator 
a polynomial in the derivatives of $A$ up to order $2$ and $\phi$ 
to order $4$. 
In summary, the results of this paper show that the 
computational apparatus and formulae for higher order terms 
given in~\cite{raichev-wilson-higher-order} hold in the case of 
complex phase functions integrated over stratified spaces.

\bibliographystyle{alpha}

\bibliography{RP}

\end{document}